\documentclass[12pt,a4paper]{amsart}
\usepackage{amsmath}
\usepackage{amssymb}
\usepackage{amsfonts}
\usepackage{mathrsfs}
\usepackage{graphicx}
\usepackage{epsfig}
\usepackage{etex}
\usepackage{multicol, color}
\usepackage[top=1.1in, bottom=1.2in,left=1.2in, right=1.2in]{geometry}
\usepackage[T1]{fontenc}
\numberwithin{equation}{section}       % Number formulas within sections

\theoremstyle{plain}

\newtheorem{prop}{Proposition}[section]
\newtheorem{theorem}{Theorem}[section]

\newtheorem{coro}[prop]{Corollary}
\newtheorem{lemma}[prop]{Lemma}

\usepackage{hyperref}
\catcode`@=11 \catcode`!=11

\expandafter\ifx\csname fiverm\endcsname\relax
  \let\fiverm\fivrm
\fi
  
\let\!latexendpicture=\endpicture 
\let\!latexframe=\frame
\let\!latexlinethickness=\linethickness
\let\!latexmultiput=\multiput
\let\!latexput=\put
 
\def\@picture(#1,#2)(#3,#4){%
  \@picht #2\unitlength
  \setbox\@picbox\hbox to #1\unitlength\bgroup 
  \let\endpicture=\!latexendpicture
  \let\frame=\!latexframe
  \let\linethickness=\!latexlinethickness
  \let\multiput=\!latexmultiput
  \let\put=\!latexput
  \hskip -#3\unitlength \lower #4\unitlength \hbox\bgroup}

\catcode`@=12 \catcode`!=12
\input pictex.tex
\catcode`@=11 \catcode`!=11
  
\let\!pictexendpicture=\endpicture 
\let\!pictexframe=\frame
\let\!pictexlinethickness=\linethickness
\let\!pictexmultiput=\multiput
\let\!pictexput=\put

\def\beginpicture{%
  \setbox\!picbox=\hbox\bgroup%
  \let\endpicture=\!pictexendpicture
  \let\frame=\!pictexframe
  \let\linethickness=\!pictexlinethickness
  \let\multiput=\!pictexmultiput
  \let\put=\!pictexput
  \let\input=\@@input   % \@@input is LaTeX's saved version of TeX's primitive
  \!xleft=\maxdimen  
  \!xright=-\maxdimen
  \!ybot=\maxdimen
  \!ytop=-\maxdimen}

\let\frame=\!latexframe

\let\pictexframe=\!pictexframe

\let\linethickness=\!latexlinethickness
\let\pictexlinethickness=\!pictexlinethickness

\let\\=\@normalcr
\catcode`@=12 \catcode`!=12
\theoremstyle{definition}
\newtheorem{definition}[prop]{Definition}

\newtheorem{remark}[prop]{Remark}

\newtheoremstyle{citing}% name
  {3pt}%      Space above, empty = `usual value'
  {3pt}%      Space below
  {\itshape}% Body font
  {}%         Indent amount (empty = no indent, \parindent = para indent)
  {\bfseries}% Thm head font
  {.}%        Punctuation after thm head
  {.5em}%     Space after thm head: " " = normal interword space;
  {\thmnote{#3}}% Thm head spec

\theoremstyle{citing}
\DeclareMathAlphabet{\mathpzc}{OT1}{pzc}{m}{it} % Zapf Chancery math alphabet

\newcommand{\C}{\mathbb{C}}
\newcommand{\D}{\mathbb{D}}

\renewcommand{\H}{\mathbb{H}}

\newcommand{\Z}{\mathbb{Z}}

\newcommand{\teta}{\widetilde{\teta}}

\newcommand{\eps}{\varepsilon}
\renewcommand{\=}{ : = }

\newcommand{\CC}{\overline{\C}}% Riemann sphere
\def\nup{\nu}
\begin{document}
\title[]{Transversality for critical relations of families of rational maps: an elementary proof}
\author{Genadi Levin, Weixiao Shen and Sebastian van Strien}
\address{Hebrew University, Shanghai Center for Mathematical Sciences, Fudan University and Imperial College London}
\date{6 April 2017. Published in  {\em New trends in one-dimensional dynamics}, 201-220, Springer Proc. Math. Stat., 285, Springer, Cham, 2019. }

\maketitle
%\tableofcontents

\begin{center}In memory of our dear friend and colleague Welington de Melo.
\end{center}

\begin{abstract}
In this paper we will give a short and elementary proof that critical relations unfold transversally
in the space of rational maps.
\end{abstract}
\section{Introduction}
In this short paper we will give an elementary proof of some transversality properties for families of
rational maps. We will consider the space
$\textbf{Rat}^{\pmb \mu}_d$  of rational maps of degree
$d$  with precisely $\nup$ critical points of multiplicities  $(\mu_1,\mu_2,\dots,\mu_\nup)$.
In Theorem~\ref{thm:cvrational} we will show that this space of maps can be locally parametrised
by critical values.  Given $f\in \textbf{Rat}^{\pmb \mu}_d$, let $\zeta=\zeta(f)\ge 0$ be
the maximal number of critical points with pairwise disjoint infinite orbits and define $N=\nup-\zeta(f)$.
In Theorem~\ref{thm:crrational-minimal} %\ref{thm:crrational}
we will show that if $f$ is not a flexible Latt\`es map then one can organise the set of critical relations
of $f$  in the form
$$\{f^{m_k}(c_{i_k})=f^{n_k}(c_{j_k}), k=1,\dots,N\}$$
so that  the map
\begin{equation}
 \textbf{Rat}^{\pmb \mu}_d \ni g \mapsto \{\sigma(g^{m_k}(c_{i_k}(g)))-\sigma(g^{n_k}(c_{j_k}(g)))\}_{k=1}^N \label{eq1}
\end{equation}
has maximal rank for $g$ near $f$, where $\sigma$ is any M\"obius transformation with $\sigma(f^{m_k}(c_{i_k}))\not=\infty$.
Property (\ref{eq1}) obviously is a transversality condition.

In fact, the choice of critical relations is in general not unique,
but as long as the selected collection is {\em full},  as made explicit in Definition \ref{def:minimal} below,
the maximal rank property holds.

Indeed, we should emphasise that some care is required
in the choice of critical relations. For example, in the case of $f_t(z)=z^2+t$ with $t=0$,
 the derivative of $t \mapsto f_t^2(0)-f_t(0)$ vanishes at $t=0$. The correct way of
expressing transversal unfolding of the critical relation $f_t(0)=0$ in (\ref{eq1}) is by taking $m_1=1$ and $n_1=0$ in this equation,
i.e. by  asserting that the derivative $t\mapsto  f_t(0)-0$ is non-zero at $t=0$.

In the unicritical case,  transversal unfolding of critical relations in the pre-periodic case goes back to Douady-Hubbard \cite{DH1} and  Tsujii \cite{Tsu0}, see also
 \cite[Remark 5.10]{Le00}.
 Milnor-Thurston \cite{MT} and Sullivan, see \cite[Theorem VI.4.2]{MS}, proved a {\lq}topological{\rq} version of transversality.

 An abstract approach to transversality for finite type maps
was developed by A. Epstein, see \cite{Ep2, Ep3}, obtaining in Part 1 of \cite{Ep2}  transversality
within the Teichm\"uller deformation space $\mbox{Def}^B_A(f)$, and in Section 5.4 in \cite{Ep2} the loci defined by critical relations within $\mbox{Def}^B_A(f)$ is discussed.
Part 2, and in particular Section 10, of \cite{Ep2} goes into a strategy for  transferring the transversality results obtained in $\mbox{Def}^B_A(f)$ to the space of rational functions.
However, we were not able to find an explicit statement covering  Theorem~\ref{thm:crrational-minimal} or Theorem~\ref{thm:crrational}.
Nevertheless, it is likely that the strategy in \cite{Ep2} can be executed  to obtain statements similar to the ones in this paper.

Our results also hold in the setting of degenerate critical points and gives unfoldings of critical relations even
when critical points share the same critical value. For this we use that  $\textbf{Rat}^{\pmb \mu}_d$
is a manifold and that $\textbf{Rat}^{\pmb \mu}_d\ni f \mapsto (f(c_1),\dots,f(c_\nup))$ has rank $\nup$,
see Theorem~\ref{thm:cvrational}.
%in particular under the assumption that all critical points are non-degenerate and ignoring deformations of
%critical relations  of the form $g(c_i)=g(c_j)$. (The approach in \cite{Ep2} seems to require that critical orbits sharing the same
%critical value are identified, having the downside that in that approach such critical relations cannot be unfolded transversally.)
%
%
%

In this short and self-contained paper we prove transversality following the approach developed by Levin in~\cite{Le0,Le1,Le,Le3}, see also \cite{LSY}.
The starting point of this paper are calculations from \cite{Le1,Le3} which show that
 if the transversality property (\ref{eq1}) fails at $g=f$, then one can construct a non-zero integrable meromorphic quadratic differential that is invariant under push-forward by $f$, which in turn implies that $f$ is a flexible Latt\`es example. Indeed, the main Theorem~\ref{thm:crrational} can be proved as in \cite{Le3}, see Remark~\ref{rem:levin}, although we shall provide a more direct and shorter proof in this paper.

The idea of using quadratic differentials appeared first in Thurston's characterization of post-critically finite branched covering of the 2-sphere \cite{DH2}. It has been used
%to obtain transversality results
in for example \cite{Ep1,LSY} and this was also used in \cite{Tsu} to obtain a similar statement to ours for the quadratic case.

%An abstract approach to transversality for finite type maps
%was developed by A. Epstein, see \cite{Ep2, Ep3}, obtaining in Part 1 of \cite{Ep2}  transversality
%within the Teichm\"uller deformation space $\mbox{Def}^B_A(f)$, and in Section 5.4 in \cite{Ep2} the loci defined by critical relations within $\mbox{Def}^B_A(f)$ is discussed.
%Part 2, and in particular Section 10, of \cite{Ep2} is dedicated to transferring the transversality results obtained in $\mbox{Def}^B_A(f)$ to the space of rational functions.
%However, we were not able to find an explicit statement covering  Theorem~\ref{thm:crrational-minimal} or Theorem~\ref{thm:crrational}.

Theorem~\ref{thm:crrational-minimal}   was proved previously for the case that critical points are non-degenerate and eventually mapped
into repelling periodic orbits, but never into a critical point, see \cite{Str,BE} and also \cite[Theorem 4.8]{FG}.

Transversality also holds in other settings. For example, if each critical point is mapped into a hyperbolic set, see \cite{Str},
%\footnote{Adam Epstein and Xavier Buff
%pointed out that in the conclusion of the Main Theorem in \cite{Str} $n$ should be $N$,
%and on  \cite[page 1, line -3]{Str}  {\lq}$\omega(c)$ contains no ...{\rq}
%by  {\lq}$P(f)$, defined below, contains no ...{\rq}.}
%%{\lq}moreover parabolic or critical points{\rq} should be replaced by {\lq}for any critical point $c$ in the Julia set, iterates of $c$ and $\omega(c)$ contain no critical point{\rq}.}
%%the assumption should be added that if  $\omega(c)$ is  not a periodic orbit, then $\omega(c)$ is not accumulated by iterates of
%%some recurrent critical point to ensure that $c$ is mapped into a hyperbolic set.}
when a summability condition holds along the orbit of  critical values, see~\cite{Le3,Astor}, for the unfolding of multipliers of periodic orbits,
see \cite{Le1,Ep2} and for a large class of interval maps, see \cite{LSS}.

As mentioned, the aim of this  paper is to present a proof of transversality for rational maps with critical relations in a complete and readily accessible form.

In Section~\ref{sec:poly} we discuss corresponding results for  polynomials.
%***********************

\section{Parametrising rational maps by their critical values}
\label{sec:param-rationals}

Let $\textbf{Rat}_d$ denote the collection of all rational maps of degree $d\ge 2$. This space is naturally parameterized by an open set in $P\C^{2d+1}$.

Given a non-ordered list ${\pmb \mu} = (\mu_1,\mu_2,\dots,\mu_\nup)$ with
$\sum_{i=1}^\nup \mu_i = 2d-2$, we say a rational map $f \in \textbf{Rat}_d$ is in the class $\textbf{Rat}^{\pmb \mu}_d$
if $f$ has precisely $\nup$ distinct critical points $c_1, c_2, \dots , c_\nup$ with multiplicities $\mu_1,\mu_2,\dots,\mu_\nup$ respectively.
Taking ${\bf 1}=(1,\dots,1)$, $\textbf{Rat}^{\bf 1}_d$ corresponds to the space of rational maps
with $2d-2$ non-degenerate critical points.

%\begin{remark}\label{rem:thm21}
Rational maps are not fully determined by their critical values (not even on small open subsets $W
\subset \textbf{Rat}^{\pmb \mu}_d$), because one can precompose
a rational map by a M\"obius transformation without changing its critical values.
%However, by applying the previous theorem, it is easy toM\"obius
%locally normalise maps $f$ so that
%the critical values locally parametrise all maps  which satisfy
%this normalisation from a neighbourhood of $f$.  %, see Corollary \ref{cor:normalise}.
However one can find a neighbourhood $W$ of $f$ and a normalisation
(based on precompositions with M\"obius transformations)
so that critical values parametrise all maps in $W$ satisfying this normalisation:
%\end{remark}

\begin{theorem}\label{thm:cvrational}
For each ${\pmb \mu}$,   $\textbf{Rat}^{\pmb \mu}_d$ is an embedded submanifold of dimension $\nup+3$
of $\textbf{Rat}_d$ and the functions defined by the critical values form a {\bf partial} holomorphic local coordinate system, i.e.
the map $\textbf{Rat}^{\pmb \mu}_d\ni f \mapsto (f(c_1),\dots,f(c_\nup))$ has rank $\nup$ and
can be completed by $3$ other coordinates to be a holomorphic coordinate system.
\end{theorem}
\begin{remark} Theorem~\ref{thm:cvrational} is not new. Similar statements are proved e.g. in \cite{EL}, \cite{Le} (see also \cite{Er}) using the Measurable Riemann Mapping Theorem with dependence on parameters; the idea of those proofs goes back probably to \cite{T}. Our proof borrows an idea of Douady and Sentenac \cite[Appendix A]{MTr}, and is short and elementary.
The case $\mu_\nup=d-1$ corresponds to the polynomial case, which in some real cases was dealt with in \cite[p120]{MS} and \cite{MTr}, see also
\cite{EG}. \end{remark}

Theorem~\ref{thm:cvrational} follows from Proposition~\ref{prop:cvrational} below.
Assume without loss of generality (by post and pre composing $f$ by M\"obius transformations if necessary) that the critical points and the critical values
avoid the point at $\infty$. Then for each $i = 1,2,\dots,\nup$,
$$f'(c_i) = f''(c_i) = \dots = f^{(\mu_i)}(c_i) = 0, f^{(\mu_i+1)}(c_i)\ne  0.$$
Applying the Implicit Function Theorem to the maps $(g,\zeta_i)\mapsto g^{(\mu_i)}(\zeta_i)$ for $(g,\zeta_i)$ near $(f,c_i)$,
gives that there exists a neighborhood $W$ of $f$ in $\textbf{Rat}_d$ and uniquely defined functions
$\zeta_i \colon W \to \C$ which are holomorphic  such that $\zeta_i(f) = c_i$ and $g^{(\mu_i)}(\zeta_i(g)) = 0, g^{(\mu_i+1)}(\zeta_i(g))\ne 0$ for each $g \in W$.
Replacing $W$ by a smaller neighborhood, for each $g\in W$ the equation $g'(\zeta)=0$
has $\mu_i$ solutions $\zeta$ (counting multiplicity) near $c_i$. It follows that
 for any $g\in W\cap \textbf{Rat}_d^{\pmb \mu}$, \,\,  $g'(\zeta)=0$
has a unique solution near $c_i$ (with multiplicity $\mu_i$); hence
  $\zeta_i(g)$ is the only critical point of $g\in W\cap \textbf{Rat}_d^{\pmb \mu}$ near $c_i$ and it has multiplicity $\mu_i$.

For $g\in W$, write
$$\zeta_i^0(g) = g(\zeta_i(g)), \zeta_i^1(g) = g'(\zeta_i(g)), \zeta_i^2(g) = g''(\zeta_i(g)), \dots.$$
Thus $\zeta_i(g)$ is a critical point of $g$ with multiplicity $\mu_i$ if and only if
$\zeta_i^j(g)=0$ for all $1\le j\le \mu_i-1$ (note that
$g^{(\mu_i)}(\zeta_i(g)) = 0, g^{(\mu_i+1)}(\zeta_i(g))\ne 0$ holds automatically for all $g\in W$).
Define $G\colon W\to \C^{2d-2}$ by
$$g  \to (\zeta_1^0(g), \zeta_1^1(g), \dots, \zeta_1^{(\mu_1-1)}(g), \dots, \zeta_\nup^0(g), \zeta_\nup^1(g), \dots , \zeta_\nup^{(\mu_\nup-1)}(g)).$$
Since $W$ has dimension $2d+1$, Theorem~\ref{thm:cvrational} follows immediately from:
\begin{prop}\label{prop:cvrational}
For each rational map $f$ as above, the Jacobian of $G$ has rank $2d -2$ at $g = f$.
\end{prop}

This proposition also immediately implies:

\begin{coro} Assume that all critical points of $f$ are non-degenerate. Then there exists a neighbourhood $W$ of $f$ in $\textbf{Rat}_d$ so that
the critical points $c_1(g),\dots,c_{2d-2}(g)$ of $g$ depend holomorphically on $g\in W$ and
the Jacobian of the map $$g\mapsto (g(c_1(g)), g(c_2(g)), \cdots, g(c_{2d-2}(g)))$$ has maximal rank at every $g\in W$.
\end{coro}

\subsection{Proof of Proposition~\ref{prop:cvrational}}

\begin{proof}[Proof of Proposition~\ref{prop:cvrational}]  Arguing by contradiction, assume that the assertion of the proposition is false. Then there exist complex numbers $A_i^j$, $1\le i\le \nup$, $0\le j<\mu_i$, not all equal to zero, such that all partial derivatives of the map
$$\textbf{G}(g)=\sum_{i=0}^\nup\sum_{j=0}^{\mu_i-1} A_i^j \zeta_i^{(j)}(g)$$
are equal to zero at $g=f$. This means that for any holomorphic curve $f_t$ in $\textbf{Rat}_d$,
 passing through $f$ at $t=0$, the map $G(t)=\textbf{G}(f_t)$ satisfies $G'(0)=0$.
Let us write
$$f_t(z)=\frac{\sum_{k=0}^d a_k(t) z^k}{\sum_{k=0}^d b_k(t) z^k}=:\frac{P_t(z)}{Q_t(z)},$$
where $a_k(t), b_k(t)$ are holomorphic in a neighborhood of $0$ and $P_0$ and $Q_0$ are co-prime polynomials.
For $1\le i \le \nup$, $j=0,\dots,\mu_i-1$ define $v_{i,j}(t)=\zeta_i^{(j)}(f_t)$. Then
$$v_{i,j}'(0)=\left.\left(\frac{\sum_{k=0}^d a_k'(0) z^k Q_0(z)-\sum_{l=0}^d b_l'(0)z^lP_0(z)}{Q_0(z)^2}\right)^{(j)}\right|_{z=c_i},$$
where we use $f^{(j+1)}(\zeta_i(f))=f^{(j+1)}(c_i)=0$.
So
\begin{multline}\label{eqn:G'0}
0=G'(0)=\sum_{i,j} A_i^j v_{i,j}'(0)\\
=\sum_{i,j}A_i^j \left.\left(\frac{\sum_{k=0}^d a_k'(0) z^k Q_0(z)-\sum_{l=0}^d b_l'(0)z^lP_0(z)}{Q_0(z)^2}\right)^{(j)}\right|_{z=c_i}.
\end{multline}

We claim that for any polynomial $T$, we have
\begin{equation}\label{eqn:T}
\sum_{i,j} A_i^j \left.\left(\frac{T(z)}{Q_0(z)^2}\right)^{(j)}\right|_{z=c_i}=0.
\end{equation}
To see this, first notice that since $T_0(z)=\prod_{i=1}^\nup (z-c_i)^{\mu_i}$ has a zero at $z=c_i$ of multiplicity $\mu_i$, the equation (\ref{eqn:T}) holds for $T=T_0$ and  $T=T_0 U$, where $U$ is an arbitrary polynomial.
 Since $\deg(T_0)=2d-2$ and any polynomial can be written as $T_0 U+T$ where $\deg(T)<2d-2$ it therefore suffices to prove (\ref{eqn:T}) in the case that $\deg(T)<2d-2$. %Therefore it suffices to prove (\ref{eqn:T}) in the case that $\deg(T)<2d-2$.
For such a polynomial $T$, we can find  polynomials $R, S$ of degree at most $d-1$ such that
$T=RQ_0-SP_0$, since $P_0$ and $Q_0$ are coprime and one of them has degree $d$. Choosing $a_k, b_l$ suitably such that $R(z)=\sum_k a_k'(0) z^k$ and $S=\sum_l b_l'(0) z^l$ and applying (\ref{eqn:G'0}), we obtain (\ref{eqn:T}).
%$$\sum_{i,j} A_i^j \left.\left(\frac{R(z) Q_0(z)-S(z) P_0(z)}{Q_0(z)^2}\right)^{(j)}\right|_{z=c_i}=0,$$
%since $a_k'(0), b_l'(0)$ are arbitrary. %we have
%Since $P_0$ and $Q_0$ are co-prime and $\max(\deg(P_0),\deg(Q_0))=d$, for any polynomial $T(z)$ of degree at most $2d$, we can find $R(z), S(z)$ of degree at most $d$ such that %$R(z)Q_0(z)-S(z) P_0(z)=T(z)$. Therefore

%For each $j$, $(T^{(j)}(z)-(T(z)/Q_0(z)^2)^{(j)})|_{z=c_i}$ can be written as a linear combination of $Q_0(c_1)^2(T(z)/Q_0(z)^2)^{j'}|_{z=c_i}$, $0\le j'<j$. It follows that there exist complex numbers $B_i^j$, not all equal to zero, such that
%$$\sum_{i,j} B_i^j T^{(j)}(c_i)=0$$
%holds for all polynomials $T$ of degree at most $2d$.
We shall now deduce from this equation that $A_i^j=0$ for all $i,j$ and thus obtain a contradiction.
Indeed, (\ref{eqn:T}) implies that for any polynomial $V$, we have
$$\sum_{i,j} A_i^j V^{(j)}(c_i)=0.$$
Fix $1\le i_0\le \nup$, $1\le j_0<\mu_{i_0}$, take
$$V(z)= \prod_{i\not= i_0} (z-c_i)^{\mu_i} (z-c_{i_0})^{j_0}.$$
Then $V^{(j_0)}(c_{i_0})\not=0$ and $V^{(j)}(c_i)=0$ for any other $(i,j)$. Therefore $A_{i_0}^{j_0}=0$.  The proof is completed.
\end{proof}

\section{Transversality results for rational maps}
Throughout this section we again consider a map $f$ in the space $\textbf{Rat}^{\pmb \mu}_d$ of rational maps of degree $d$, with
$\nu$ distinct critical points $c_1, c_2, \dots , c_\nup$ with multiplicities ${\pmb \mu}=(\mu_1,\mu_2,\dots,\mu_\nup)$ where $\sum_{i=1}^\nup \mu_i = 2d-2$. For $g$ in a small neighborhood of $f$ in $\textbf{Rat}^{\pmb \mu}_d$, the critical points $c_1(g), c_2(g), \dots, c_\nup(g)$ depends holomorphically on $g$.

We are interested in the smoothness of sets defined by a set of critical relations of the form $g^m(c_i(g))=g^n(c_j(g))$. A particular case of our main result in this direction is the following:
\begin{theorem}\label{thm:rationalsimple} Let $f\in \textbf{Rat}^{\,\pmb \mu}_d$ and assume that there exists $1\le i, j\le \nup$ and $m>0$
so that $f^m(c_i)=c_j$. % and $f^k(c_i)$ is not a critical point for $0<k<m$.  
Then the equation
$$g^m(c_i(g))=c_j(g)$$
defines an embedded
submanifold of $\textbf{Rat}^{\pmb \mu}_d$ of codimension one near $f$.
% unless $f$ is a flexible Latt\'es example. (The last possibility is automatically rules out if $n=0$.)
\end{theorem}

In order to state a more general result, we have to prepare some terminology. Let us say that a quadruple $(i,j;m,n)$ is a {\em (candidate) critical relation} if $1\le i, j\le \nu$, and $m,n$ are non-negative integer with $m+n>0$.
We say that this critical relation is {\em realized by $f$} if $f^m(c_i(f))=f^n(c_j(f))$.

%Assume without loss of generality that the critical orbits of $f$ avoid the point at $\infty$,
%see Remark~\ref{sigma}. \marginpar{NEW}
%For each finite {\em collection} $\mathcal{F}=\{(i_k,j_k; m_k,n_k): 1\le k\le N\}$ of critical relations realized by $f$, consider the map
%\begin{equation}\label{eqn:crmulti}
%\mathcal{R}_{\mathcal{F}}: g\mapsto \left( g^{m_k}(c_{i_k}(g))-g^{n_k}(c_{j_k}(g))\right)_{k=1}^N
%\end{equation}
%which is holomorphic in a neighborhood of $f$ in $\textbf{Rat}^{\pmb \mu}_d$.

Given $f$, let $\zeta=\zeta(f)\ge 0$ be
the maximal number of critical points with pairwise disjoint infinite orbits.  Note that this number is well-defined,
but that one {\em cannot} say {\em which} critical points are {\lq}free{\rq}.
For example, if $f$ has three distinct critical points $c_1,c_2,c_3$, so that the forward orbits of
$f(c_1)=f(c_2)$ and $c_3$  are disjoint and infinite, then $\zeta(f)=2$; of course one could consider $c_1,c_3$ as
the free critical points of $f$, but equally well also $c_2,c_3$.

In this section we will show

\begin{theorem}\label{thm:crrational-minimal}
Assume $f\in \textbf{Rat}_d^{\pmb \mu}$ is not a flexible Latt\`es map. Then  there exists a set
$$\mathcal F=\{(i_k, j_k; m_k, n_k), k=1,\dots,N\}\mbox{ with }N=\nu-\zeta(f)$$
of critical relations $f^{m_k}(c_{i_k})=f^{n_k}(c_{j_k})$ which are realised by $f$,
such that the Jacobian of the map
\begin{equation}\mathcal{R}_{\mathcal{F}}^\sigma \colon g\mapsto (\sigma(g^{m_k}(c_{i_k}(g)))-\sigma(g^{n_k}(c_{j_k}(g))))_{k=1}^N
\label{eqn:crmulti}\end{equation}
at $g=f$ has rank $N$, whenever $\sigma$ is a M\"obius transformation  for which $\sigma(f^{m_k}(c_{i_k}))\in \C$, $k=1,\dots,N$.
%the Jacobian of the map
%$g\mapsto (g^{m_k}(c_{i_k}(g))-g^{n_k}(c_{j_k}(g)))_{k=1}^N$
%at $g=f$ has rank $N$.
%Assume $f\in \textbf{Rat}_d^{\pmb \mu}$ is not a flexible Latt\`es map and the critical orbits avoid the point at $\infty$. Then  there exists a set
%$$\mathcal F=\{(i_k, j_k; m_k, n_k), k=1,\dots,N\}\mbox{ with }N=\nu-\zeta(f)$$
%of critical relations $f^{m_k}(c_{i_k})=f^{n_k}(c_{j_k})$ which are realised by $f$, so that the %Jacobian of the map $\mathcal{R}_{\mathcal{F}}$ at $g=f$ has rank $N$.
%\begin{equation}\label{eqn:crmulti}
% \textbf{Rat}^{\pmb \mu}_d \ni g \mapsto \{g^{m_k}(c_{i_k}(g))-g^{n_k}(c_{j_k}(g))\}_{k=1}^N
%\end{equation}
%has rank $N$ for all $g$ close to $f$.
\end{theorem}

\begin{remark}\label{sigma}
The assumption that $\sigma(f^{m_k}(c_{i_k}))\in \C$ is made to ensure that (\ref{eqn:crmulti}) is holomorphic near $f$.
The kernel of the Jacobian of $\mathcal{R}_{\mathcal{F}}^\sigma$ at $f$, hence its rank, does not depend on $\sigma$, as long as $\sigma(f^{m_k}(c_{i_k}))\not=\infty$ for all $k=1,\dots,N$. Indeed, a tangent vector of $\textbf{Rat}_d^{\pmb \mu}$ at $f$ belongs to the kernel if and only if it has the same image under the tangent map of the maps $g\mapsto (g^{m_k}(c_{i_k}(g)))_{k=1}^N$ and $g\mapsto (g^{n_k}(c_{j_k}(g)))_{k=1}^N$ at $g=f$ (both are holomorphic maps from a neighborhood of $f$ in $\textbf{Rat}_d^{\pmb \mu}$ into $\CC^N)$.

In particular, to prove Theorem~\ref{thm:crrational-minimal}, we can and will assume that the critical obits of $f$ avoid $\infty$  and only prove that $\mathcal{R}_{\mathcal{F}}=\mathcal{R}_{\mathcal{F}}^{id}$ has rank $N$ at $g=f$. Indeed, we can always choose $z_0$ (arbitrarily close to $\infty$) which avoids the critical orbits of $f$. Put $\sigma(z)=z_0z/(z_0-z)$ and $\tilde{f}=\sigma \circ f\circ \sigma^{-1}$. Then $\infty$ avoids the critical orbits of $\tilde{f}$. Since
$\mathcal{R}_{\mathcal{F}}^\sigma(g)=\mathcal{R}_{\mathcal{F}}^{id} (\sigma \circ g\circ \sigma^{-1}),$
once we prove that the Jacobian of $g\mapsto\mathcal{R}_{\mathcal{F}}^{id}(g)$ has rank $N$ at $g=\tilde{f}$, it follows that the Jacobian of $\mathcal{R}_{\mathcal{F}}^\sigma$ has rank $N$ at $g=f$.
\end{remark}

\begin{remark}\label{rem:minimalfullThm32}
There are several ways of assigning a set of critical relations $\mathcal F$ to $f$.
%One construction for $\mathcal F$ is given in Lemma~\ref{lem:comb} below.
As we will prove in Subsection~\ref{subsec:pf32}, for {\em any} set of critical relations which is {\em full} in the sense of Definition~\ref{def:minimal},
Theorem~\ref{thm:crrational-minimal} holds. %\marginpar{minimal only required}
\end{remark}

\begin{remark} A {\em flexible Latt\`es map} is by definition a rational map that is conformally conjugate to a map of the form $L/ \!\! \sim\colon T  /\!\! \sim \to T / \!\! \sim$, where $T=\C/ (\Z \oplus \gamma \Z)$, $\gamma\in \H$ (where $\H$ is the upper-half plane), $\sim$ is the equivalence relation on $\C$ defined  by $z\sim -z$
 and $L\colon \C\to \C$ is of the form $L(z)=az+b$ with $a\in \Z$ and $2b\in \Z \oplus \gamma \Z$, see \cite{Mil}.
Such maps can be of two types: either each critical point is mapped in two iterates into a repelling fixed point
 or in one iterate into a repelling periodic point of period two, see  \cite{Mil}.
  \end{remark}

\begin{remark}
Theorem \ref{thm:crrational-minimal} and the implicit function theorem, imply that manifolds defined by critical relations corresponding to disjoint
subsets $\mathcal F',\mathcal F''$ of $\mathcal F$ are smooth
and transversal to one another.
\end{remark}

For completeness we prove the following corollary of Theorem~\ref{thm:crrational-minimal}:

\begin{coro}
If each critical point $c_i$ is eventually mapped to a repelling periodic point $p_i$ with $f^{m_i}(c_i)=p_i$ and
 $f^j(c_i)\notin \{c_1,\dots,c_\nu\}$ for all $j=1,\dots,m_i$ then the Jacobian of
\begin{equation}\label{eqn:crmultip}
 \textbf{Rat}^{\pmb \mu}_d \ni g \mapsto \{\sigma(g^{m_i}(c_i(g)))-\sigma(p_i(g))\}_{i=1}^\nu
\end{equation}
has maximal rank at $g=f$, where $\sigma$ is a M\"obius transformation with $\sigma(p_i)\not=\infty$ for all $i$.
\end{coro}
\begin{proof} For the same reason as explained in Remark~\ref{sigma}, we only need to consider the case where $\infty$ avoids the critical orbits and $\sigma=id$. Let $\mathcal{R}$ denote the map in (\ref{eqn:crmultip}).
The corollary follows from the following claim by Theorem~\ref{thm:crrational-minimal}.

{\bf Claim.} If $f_t$ is a holomorphic curve in $\textbf{Rat}^{\pmb \mu}_d$ passing through $f$ at $t=0$ which represents a vector in the kernel of $D_f\mathcal{R}$, then for any critical relations $(i,j;m,n)$ realized by $f$, we have $$f_t^m(c_i(f_t))-f_t^n(c_j(f_t))=o(t)\text{ as } t\to 0.$$

Indeed, the claim implies that the kernel of $D_f\mathcal{R}$ is contained in the kernel of $D_f\mathcal{R}_{\mathcal F}$ for any finite collection $\mathcal{F}$ of critical relations.  By Theorem~\ref{thm:crrational-minimal}, we can choose $\mathcal{F}$ such that $D_f\mathcal{R}_{\mathcal{F}}$ has maximal rank. Thus $D_f\mathcal{R}$ has maximal rank.

Let us prove the claim. Choose $k$ large enough such that $m+k\ge m_i$ and $n+k\ge m_j$.
Since $f^{m+k-m_i}(p_i)=f^{n+k-m_j}(p_j)$, the periodic points $p_i$ and $p_j$ have the same period, denoted by $s$.
Moreover, if $p_i(t)$ (resp. $p_j(t)$) denotes the repelling periodic point of $f_t^s$ near $p_i$ (resp. $p_j$), then
$$f_t^{m+k-m_i}(p_i(t))=f_t^{n+k-m_j}(p_j(t)).$$
Since $f_t^{m_i}(c_i(f_t))-p_i(t)=o(t)$ as $t\to 0$, we have
$$f_t^{m+k}(c_i(f_t))-f_t^{m+k-m_i}p_i(t)=o(t)\text{ as }t\to 0.$$ Similarly, we have
$$f_t^{n+k}(c_j(f_t))-f_t^{n+k-m_j}(p_j(t))=o(t)\text{ as }t\to 0.$$ Therefore,
$$f_t^{m+k}(c_i(f_t))-f_t^{n+k}(c_j(f_t))=o(t) \text{ as } t\to 0.$$
Since $f^{m+k'}(c_i)$ is not critical for each $0\le k'<k$, it follows that
$f_t^{m}(c_i(f_t))-f_t^n(c_j(f_t))=o(t)$ as $t\to 0$.
\end{proof}

\subsection{How to associate critical relations to a rational map}\label{subsec:howto}
There are several ways to record the (infinitely many) critical relations of a rational map.
In this subsection we will show how one can associate these
in an efficient way so that in particular no critical relation is counted twice.

As above, let $c_1,c_2,\ldots, c_\nup$ be the critical points of a rational map in the class $\textbf{Rat}_d^{\pmb \mu}$.
For an arbitrary collection $\mathcal{F}$ of critical relations realized by $f$, let $\sim_{\mathcal{F}}$ denote the
{\em smallest} equivalence relation in the set
$\Sigma:=\{(i,m): 1\le i\le \nup, m\ge 0\}$ such that $(i,m+k)\sim_{\mathcal{F}} (j,n+k)$ for each $(i,j;m,n)\in\mathcal{F}$ and each $k\ge 0$.

So  $\sim_{\mathcal{F}}$ defines the set of critical relations that can be {\lq}read off{\rq}  from $\mathcal F$.
So for example, if $\nu=4$ and $\mathcal F= \{(1,2; 1,1), (1,3;1,1)\}$ then $(i,1+k)\sim_{\mathcal{F}} (j,1+k)$ for all $i,j\in \{1,2,3\}$ and all $k\ge 0$, but $(i,m)\not  \sim_{\mathcal{F}} (4,n)$ for $i\in \{1,2,3\}$ and all $m,n\ge 0$.

%, one can associate a directed graph $(V_f,E_f)$
\iffalse
whose vertices $V_f$ consists
of forward iterates of the critical points, i.e. $V_f=\{f^n(c_i); 1\le i\le \nup, n\ge 0\}$  and
with an edge $a\to b$ between $a,b\in V_f$  if and only if $b=f(a)$.

where $V_f=\{(i,m): 1\le i\le \nup, n\ge 0\}$ and with an edge $(i,m)\to (j,n)$ if and only if $f^{m+1}(c_i)=f^n(c_j)$.
Any edge of the form $(i,n)\to (i,n+1)$ is called {\em trivial}.

Similarly, a set of critical relations $\mathcal{F}$ also defines a directed graph $(V_f,E_{\mathcal F})$
where $E_{\mathcal F}$ consists of trivial edges
and, for each $(i,j;m,n)\in \mathcal{F}$, edges of the form
$(j,n-1+k)\to (i,m+k)$
(when $n\ge 1$) and $(i, m-1+k)\to (j,n+k)$ (when $m\ge 1$), $k=0,1,\ldots$.
\fi

\begin{figure}[ht]
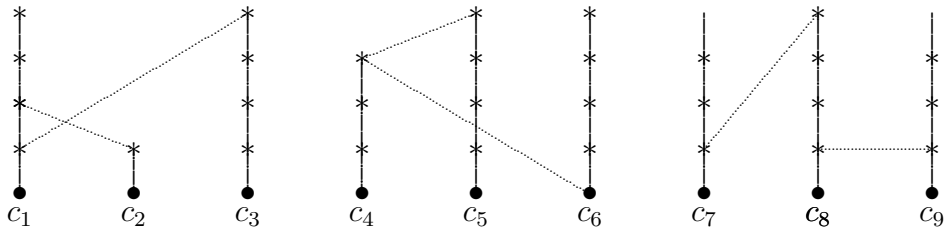

\beginpicture
\dimen0=0.15cm
\dimen1=0.06cm
%\arrow <10pt> [.2,.67] from 33 8 to 5  -26
\setcoordinatesystem units <\dimen0,\dimen1> point at 0 0
\setplotarea x from 15 to 80, y from 0 to 50
\setlinear
\plot  10 0  10 40 /
\put {$c_1$} at 10 -5
\put {$\bullet$} at 10 0
\put {$\ast$} at 10 10
\put {$\ast$} at 10 20
\put {$\ast$} at 10 20
\put {$\ast$} at 10 30
\put {$\ast$} at 10 40
%\put {$\ast$} at 10 50
\plot  20 0  20 10 /
\put {$c_2$} at 20 -5
\put {$\bullet$} at 20 0
\put {$\ast$} at 20 10
%\put {$\ast$} at 20 20
\plot  30 0  30 40 /
\put {$c_3$} at 30 -5
\put {$\bullet$} at 30 0
\put {$\ast$} at 30 10
\put {$\ast$} at 30 20
\put {$\ast$} at 30 30
\put {$\ast$} at 30 40
\put {$c_4$} at 40 -5
\plot  40 0  40 30 /
\put {$\bullet$} at 40 0
\put {$\ast$} at 40 10
\put {$\ast$} at 40 20
\put {$\ast$} at 40 30
%\put {$\ast$} at 40 40
\plot  50 0  50 40 /
\put {$c_5$} at 50 -5
\put {$\bullet$} at 50 0
\put {$\ast$} at 50 10
\put {$\ast$} at 50 20
\put {$\ast$} at 50 30
\put {$\ast$} at 50 40
\plot  60 0  60 40 /
\put {$c_6$} at 60 -5
\put {$\bullet$} at 60 0
\put {$\ast$} at 60 10
\put {$\ast$} at 60 20
\put {$\ast$} at 60 30
\put {$\ast$} at 60 40
\setdots <0.5mm>
\plot 20 10 10 20 /
\plot 30 40 10 10 /
\plot 40 30 60 00 /
\plot 50 40 40 30 /
\setsolid
\put {$\bullet$} at 70 0
\put {$\ast$} at 70 10
\put {$\ast$} at 70 20
\put {$\ast$} at 70 30
\put {$c_7$} at 70 -5
\put {$\bullet$} at 80 0
\put {$\ast$} at 80 10
\put {$\ast$} at 80 20
\put {$\ast$} at 80 30
\put {$\ast$} at 80 40
\put {$\bullet$} at 90 0
\put {$\ast$} at 90 10
\put {$\ast$} at 90 20
\put {$\ast$} at 90 30
\put {$c_8$} at 80 -5
\put {$c_8$} at 80 -5
\put {$c_9$} at 90 -5
\plot 70 0 70 40 /
\plot 80 0 80 40 /
\plot 90 0 90 40 /
\setdots <0.5mm>
\plot 80 10 90 10 /
\plot 70 10 80 40 /
%\arrow <10pt> [.2,.67] from 20 10 to 10  20
%\put {$\ast$} at 30 20
\endpicture
\caption{\label{fig:firstreturn}
\small The orbit diagram of a map with $\zeta(f)=3$ and  $f^2(c_1)=f(c_2), f(c_1)=f^4(c_3),f^3(c_4)=f^4(c_5)=c_6$, $f(c_7)=f^4(c_8)$, $f(c_8)=f(c_9)$.
%Following Definiton~\ref{def:lemcomb} and Lemma~\ref{lem:comb}, the critical relations in this map are recorded as $(2,1,1,2)$, $(3,1,4,1)$, $(4,6,3,0)$, $(5,4,4,3)$.}
Each of the collections $\{(2,1;1,2),(3,1;4,1),(4,6;3,0),(5,4;4,3),(8,7;4,1),(9,8;1,1)\}$,
 $\{(2,1;1,2),(3,1;4,1), (4,6;3,0),(5,6;4,0),(9,7;4,1),(9,8;1,1)\}$ and
 $\{(2,1;2,3),(3,1;5,2), (4,6;3,0),(5,6;4,0),(9,7;4,1),(9,8;2,2)\}$ is minimally full.}
\end{figure}

Roughly speaking, we say that a collection $\mathcal{F}$ of critical relations
is {\em full} if it  {\lq}essentially{\rq} explains all critical relations of
$f$ and  $\mathcal{F}$ is {\em minimally full} if it does not contain redundant
critical relations.  More precisely,
%
%We say that a critical relation $(i,j;m,n)$ realized by $g$ is  {\em minimal}, if all the points
%$g^k(c_i),g^l(c_j)$, $0\le k\le m-1,0\le l \le n-1$ are distinct (if $m=0$ then we merely assume
%that $g^l(c_j)$, $0\le l \le n-1$ are disjoint; and similarly for $n=0$).
\begin{definition}\label{def:minimal} We say that
a collection $\mathcal{F}$ of critical relations realized by $f$ is {\em full} if for
any critical relation $(i,j;m,n)$ realised by $f$,
i.e. whenever $f^m(c_i)=f^n(c_j)$,
there exists $k\ge 0$ such that  $(i, m+k)\sim_{\mathcal F} (j, n+k)$
and such that $f^{m+k'}(c_i)=f^{n+k'}(c_j)\not\in \{c_1, \dots, c_\nup\}$ for each $0\le k'<k$.

Note that any full collection contains at least $\nup-\zeta(f)$ relations.
% candidate critical relation $(i,j;m,n)$ realized by $f$, we have $(i,m)\sim_{\mathcal{F}} (j,n).$
A full collection $\mathcal{F}$ is called {\em minimally full} if
$\#\mathcal{F}=\nup-\zeta(f)$.
%it consists of $\nup-\zeta(f)$ critical relations.
%\begin{itemize}
%\item {\em minimal} if there exists no proper subset $\mathcal{F}_0$ of $\mathcal{F}$
%the equivalence relations $\sim_{\mathcal F}$ and $\sim_{{\mathcal F}_0}$ coincide; and
%\item {\em full}

If $\mathcal F$ is minimally full then in particular there  exists no $1\le i_1, i_2, \cdots, i_k\le \nup$, $k\ge 2$, such that
\begin{equation}
(i_1,i_2;1,1), (i_2, i_3; 1,1),\cdots, (i_k,i_1;1,1)\in \mathcal{F}.\label{eq:min11}
\end{equation}
We refer to the last property as the {\em non-cyclic condition}.
\end{definition}

So if $f$ has critical points  $c_1, \dots, c_4$ with critical relations $f^k(c_1)=f^k(c_2)=f^k(c_3)$ for all $k\ge 1$ and there are no other critical relations,
 then  $\zeta(f)=2$ and
$$\mathcal F_1= \{(1,2; 1,1), (1,3;1,1)\} \mbox{ but also }\mathcal F_2= \{(1,2; 2,2), (1,3;3,3)\} $$ are minimally full collections.

Note that if $\mathcal  F$ is minimally full and $(i,j;m,n)\in \mathcal F$ then $(j,i;n,m)\notin \mathcal F$.
Later on, we will define a convenient choice for a minimally full collection $\mathcal F$, see
Definition~\ref{def:lemcomb} and in Lemma~\ref{lem:comb} we will show that such a choice can always be made.

%\iffalse
%  \marginpar{Added}
%\begin{remark}\label{rem:noncrit}
%A minimal collection of critical relations $$\mathcal F=\{(i_k, j_k; m_k, n_k), k=1,\dots,N\}$$ associated to $f$
%has by definition the property that for each $k=1,\dots,N$, both
%$f(c_{i_k}),\dots,f^{m_k-1}(c_{i_k})$ and $f(c_{j_k}),\dots,f^{n_k-1}(c_{j_k})$ are disjoint from the set of critical points.
%%So if $f^m(c_i)=c_j$ with $m>0$ minimal, then this should be included in the set of critical relations.
%\end{remark}
%\fi

\subsection{An even more general theorem} \label{subsec:32}
Let $g$ be a rational map with critical points $c_1(g),\dots,c_\nu(g)$.
Associate to each $(i,j;m,n)$ the following rational map
$$Q^g_{i, j; m,n}(z)=\sum_{r=1}^m \frac{Dg^{m-r}(g^r(c_i(g)))}{z-g^r(c_i(g))}-\sum_{s=1}^n\frac{Dg^{n-s}(g^s(c_j(g)))}{z-g^s(c_j(g))},$$
when $g^{r}(c_i(g)), g^s(c_j(g))\not=\infty$ for all $1\le r\le m, 1\le s\le n$.
(Convention: For $m=0$ or $n=0$, the corresponding sum is understood as $0$.)

Given a meromorphic quadratic differential $Q=q(z) dz^2$, define its push-forward as $f_*Q= \widehat{q}(z) dz^2$, where
$$\widehat{q}(z)=\sum_{w\in f^{-1}(z)}\frac{q (w)}{f'(w)^2}.$$
It is not difficult to check that $f_*Q$ is again a meromorphic quadratic differential. The assignment
$Q\mapsto f_* Q$ is often called the Thurston operator, see \cite{Gar},\cite{McM0}, and was used in Thurston's rigidity theorem, see \cite{DH2}.
M. Tsujii was probably the first to use quadratic differentials in the context of transversality, see \cite{Tsu0, Tsu, Mil0}, but see also
 \cite{Astor,Ep2,Ep3,Le0,Le1,Le3,Mak}.

%in his topological characterization of rational functions, see \cite{DH}.

Theorem~\ref{thm:crrational-minimal} will follow from
%the following holds:
%\begin{itemize}
%
%\end{itemize}
%if  $m\ge n\ge 0$ which are not both zero, define a rational map

%A {\em critical relation} of $f$ is an equation $f^{n_1}(c_{j_1})=f^{n_2}(c_{j_2})$, where $j_1, j_2\in \{1,2,\ldots, 2d-2\}$ (maybe equal), $n_1,n_2$ are non-nengative integers, no both %zero. The critical relation is called {\em primitive} if none of the points $f^k(c_{j_i})$, $1\le k<n_i$, $i=1,2,$, is a critical point of $f$. The {\em support} of this critical relation %is the set $\{f^k(c_{j_i}): 1\le k\le n_i, i\in \{1,2\}\}$.

%A finite collection of critical relations
%$$(A_r):  \,\,\, f^{n_1(r)}(c_{j_1(r)})=f^{n_2(r)}(c_{j_2(r)}),$$
%of $f$ is called {\em essential} if for any $r$, the support of $A_r$ is not contained in the union of supports of other critical relations in this collection.

%Let us assume that $f$ has some critical relations, formulated as follows:
%There exists $0\le r_1\le r_2\le 2d-2$, such that
%\begin{itemize}
%\item For each $1\le j\le r_1$, there exists a positive integer $q_j$ and $\nup(j)$, such that %$f^{q_j}(c_j)=c_{\nup(j)}$ and $f^k(c_j)\not\in\Crit(f)$ for $1\le k<q_j$;
%\item For $r_1<j\le r_2$, there exist positive integers $l_j<q_j$ such that %$f^{q_j}(c_j)=f^{l_j}(c_j)$ and $f^k(c_j)\not\in \Crit(f)$ for all $1\le k<q_j$.
%\end{itemize}

\begin{theorem}\label{thm:crrational}
Assume that the critical orbits of $f\in \textbf{Rat}^{\pmb \mu}_d$ avoid $\infty$. Let $\mathcal{F}$ be a finite set of critical relations $(c_{i_k}, c_{j_k}, m_k, n_k)$, $k=1,2,\ldots, N$, which are realized by $f$ and which satisfies the non-cyclic condition (\ref{eq:min11}).
%this follows from Definition~\ref{def:lemcomb}
If the Jacobian of the map
\begin{equation}\label{eqn:crmulti-bis}
 \textbf{Rat}^{\pmb \mu}_d \ni g \mapsto \{g^{m_k}(c_{i_k}(g))-g^{n_k}(c_{j_k}(g))\}_{k=1}^N
\end{equation}
at $g=f$ has rank less than $N$, then there exist complex numbers $a_1, a_2, \cdots, a_N$, such that
\begin{itemize}
\item for some $k$, $(m_k,n_k)\not=(1,1)$ and $a_k\not=0$;
\item $f_*(q(z)dz^2)=q(z) dz^2$, where
\begin{equation}\label{eqn:mainq}
q(z)=\sum_{\substack{1\le k\le N\\ (m_k,n_k)\not=(1,1)}} a_k Q^f_{{i_k}, {j_k}; m_k, n_k}(z).
\end{equation}
\end{itemize}
If in addition $f$ is not a flexible Latt\'es example, then $q(z)\equiv 0$.
%and non-negative integers $k\ge l$ such that either $i\not=j$ or $k>l\ge 0$, the following equation defines an immersed submanifold in a neighborhood of $f$ in $\textbf{Rat}_d$:
%$$g^k(c_i(g))=g^l(c_j(g)),$$
%unless $f$ is a flexible Latt\'es example.
\end{theorem}

\begin{remark} If $f$ has $2d-2$ distinct critical values, then the converse statement of the theorem also holds. Namely, for any finite set $\mathcal{F}$ as above, if $f_*(q(z)dz^2)=q(z) dz^2$, then the Jacobian of the map (\ref{eqn:crmulti}) has rank less than $N$.
\end{remark}

\begin{remark}\label{S} Take $f\in \textbf{Rat}_d^{\pmb \mu}$ and a manifold $S$ passing through $f$ of dimension $p$, that is transverse to the orbit $O(f)$ of $f$ under M\"obius conjugacies. Assume that the map defined in (\ref{eqn:crmulti})
has maximal rank. Then the restriction of this map to $S$  also has maximal rank.
This holds because the value of the  map (\ref{eqn:crmulti}) is constant on $O(f)$.
\end{remark}

\section{Theorems~\ref{thm:rationalsimple}-\ref{thm:crrational-minimal} follow from
Theorem~\ref{thm:crrational}}

\subsection{Proof of Theorem~\ref{thm:rationalsimple}}
%\begin{proof}[Proof of Theorem~\ref{thm:rationalsimple}]
If $f^m(c_i)=c_j$ then $f$ is not a Latt\'es example. We may assume without loss of generality that $\infty$ avoids the critical orbits of $f$ so that Theorem~\ref{thm:crrational} applies.
It is clear that
$$Q^f_{i, j;m, 0}(z)=\dfrac{Df^{m-1}(f(c_i))}{z-f(c_i)}+\dots + \dfrac{1}{z-f^{m}(c_i)}$$
has a pole at $c_j$, so it is not identically zero and thus the conclusion follows from the last sentence
of Theorem~\ref{thm:crrational}.
\qed
%\end{proof}

\subsection{An improved way to organise critical relations and the proof of Theorem~\ref{thm:crrational-minimal}}\label{subsec:pf32}

In general, one can associate several full collections $\mathcal F$ to $f$ each giving rise to a
map $\mathcal{R}_{\mathcal F}^\sigma$ as in (\ref{eqn:crmulti}). Let us first prove, as claimed in  Remark~\ref{rem:minimalfullThm32}, that
any full collection gives rise to the same rank:

\begin{lemma}\label{lem:samerank} For any full collections $\mathcal{F}$ and $\mathcal{F}'$ of critical relations for $f$, the Jacobian matrices of $\mathcal{R}_{\mathcal{F}}^\sigma$ and $\mathcal{R}_{\mathcal{F}'}^\sigma$ at $g=f$ have the same rank.
% The rank of derivative of
%\begin{equation} \textbf{Rat}^{\pmb \mu}_d \ni g \mapsto \{g^{m_k}(c_{i_k}(g))-g^{n_k}(c_{j_k}(g))\}_{k=1}^N
%\label{eq:jacobian} \end{equation}
%at $g=f$  does not depend on the choice of the minimal collection.
\end{lemma}
\begin{proof}
According to Remark~\ref{sigma}, we may assume the critical orbits avoid $\infty$ and $\sigma=id$.
Consider a holomorphic curve $f_t$, passing through $f$ at $t=0$.
This curve represents a vector in the kernel of $D\mathcal{R}_{\mathcal{F}}$ if and only if the derivative of
$t\mapsto f_t^{m}(c_i(f_t))-f_t^{n}(c_j(f_t))$ vanishes at $t=0$ for each $(i,j;m,n)\in\mathcal{F}$,
and therefore if and only if $t\mapsto f_t^{m}(c_i(f_t))-f_t^{n}(c_j(f_t))$ vanishes at $t=0$ for each $(i,m)\sim_{\mathcal{F}} (j,n)$
where $\sim_{\mathcal{F}}$ is the equivalence relation associated to $\mathcal{F}$ as
defined in the first paragraph of Section~\ref{subsec:howto}.

Assume that $(i,j;m,n)$ is realised by $f$.
Since $\mathcal{F}$ is full,  there exists $k\ge 0$ so that  $(i,m+k)\sim_{\mathcal F} (j,n+k)$ and so that
$Df^{k}(f^m(c_i))=Df^k(f^n(c_j))\ne 0$. So if $f_t$ represents a vector in the kernel of $D\mathcal{R}_{\mathcal{F}}$
then the derivative of $t\mapsto f_t^{m+k}(c_i(f_t)) - f_t^{n+k}(c_j(f_t))$ vanishes at $t=0$. Since
$f^m(c_i)=f^n(c_j)$ and $Df^{k}(f^m(c_i))=Df^k(f^n(c_j))\ne 0$, this implies that
the derivative of $t\mapsto f_t^{m}(c_i(f_t))-f_t^{n}(c_j(f_t))$ vanishes at $t=0$.

On the other hand,
if for each $(i,j;m,n)$ which is realised by $f$  the derivative of $t\mapsto f_t^{m}(c_i(f_t))-f_t^{n}(c_j(f_t))$ vanishes at $t=0$,
then in particular this holds for each  $(i,j;m,n)\in\mathcal{F}$ and so the
holomorphic curve $f_t$ represents a vector in the kernel of $D\mathcal{R}_{\mathcal{F}}$.

It follows that $f_t$ represents a vector in the kernel of $D\mathcal{R}_{\mathcal{F}}$ if and only
 if for each $(i,j;m,n)$ which is realised by $f$ the derivative of $t\mapsto f_t^{m}(c_i(f_t)) - f_t^{n}(c_j(f_t))$ vanishes at $t=0$.
%this holds if and only if the derivative of $t\mapsto  f_t^{m}(c_i(f_t))-f_t^{n}(c_j(f_t))$ vanishes at $t=0$, for each $(i,j;m,n)$ realized by $f$.
The last condition is independent of the choice of the full collection $\mathcal F$.
Since both $\mathcal{F}$ and $\mathcal{F}'$  are full, the rank-nullity theorem implies that the rank of the Jacobian matrices are the same.

\end{proof}

We will find it convenient to prove  Theorem~\ref{thm:crrational-minimal} for a conveniently chosen minimal collection $\mathcal F$,
namely one which satisfies the  following stronger  minimality assumption.

\begin{definition}\label{def:lemcomb}
We say that a {\em collection} $\mathcal{F}$ is {\em proper} for $f$  if it is minimally full and satisfies the following extra properties:
 \begin{enumerate}
 \item $(i,j;m,n)\in \mathcal{F}$ implies $m>0$ and either $i\ge j$ or $n=0$.  If $i=j$ then
 $m>n$.
 \item if $(i,j;m,n)\in \mathcal{F}$ then the collection of points
$f^k(c_i), k=1,\dots,m-1$
is pairwise disjoint and does not intersect $c_1,\dots,c_\nu$
nor the forward orbits of $c_1,\dots,c_{i-1}$.
\item For each $1\le i\le \nu$  there exists at most one critical relation of the form $(i,j;m,n)\in \mathcal F$.
\item For each
$1\le j\le \nu$  there exists at most one critical relation of the form $(i,j;m,0)\in \mathcal F$;
\item For each $1\le j\le \nu$  and each $n>1$ there exists at most one critical relation of the form $(i,j;1,n)\in \mathcal F$;
\item If $(i,j;m,n)\in \mathcal{F}$ with $m>1$ and $n>0$, and $(k,i;1,l)\in \mathcal{F}$ for some $k$ and $l$, then $l<m$.
\end{enumerate}
\end{definition}

\begin{remark} If the collection $\mathcal{F}$ is proper then it satisfies the non-cyclic condition \ref{eq:min11}.

\end{remark}

\begin{lemma}\label{lem:comb}
%Associated to each rational map $f\in \textbf{Rat}_d$,
There exists a proper collection of critical relations which are realised by $f$.
% so that $\mathcal F$
%is proper.
\end{lemma} \begin{proof}
% Label the critical points, so that $c_1,\dots,c_q$ are periodic and
%$c_{\nup-\zeta},\dots,c_p$ are the {\lq}free orbits{\rq}.
%For each $i=1,\dots,q$, consider the critical relation $(i,i,n_i,0)$ where $n_i$ is the period of $c_i$.
%For $i=q+1,\dots,\nup-\zeta$,  inductively define $n_i>0$ maximal  so that $c_i,f(c_i),\dots,f^{n_i-1}(c_i)$ are distinct and also distinct
%from $$c_j,,f(c_j),\dots,f^{n_j-1}(c_j) \mbox{ for } j=1,\dots,i-1.$$
%If $n_i$ is finite, then $f^{n_i}(c_i)=f^{m_i}(c_j)$ for some $m_i\ge 0$ and some $1\le j\le \nup-\zeta$
%and associate to $c_i$ the critical relation $(i,j,n_i,m_i)$. The resulting set of critical relations is realised by $f$,
%and by construction complete and minimal.
% \end{proof}
For each $i=1,\dots,\nup$,  inductively define $m_i>0$ maximal  so that $f(c_i),\dots,f^{m_i-1}(c_i)$ are distinct and also distinct
from $$\{f^k(c_j);  0\le k < m_j \mbox{ , } j=1,\dots,i-1\}\cup \{c_1,\dots,c_\nup\}.$$
(When $i=1$ we take this union to be $ \{c_1,\dots,c_\nup\}$.) If $m_i$ is finite, then there are two possibilities:

(a) $f^{m_i}(c_i)=f^{n_{j_i}}(c_{j_i})$ for some  $1\le j_i\le i$ and some finite $n_{j_i}$ with $0<n_{j_i}<m_{j_i}$.%$0<  n_{j_i}\le  m_{j_i}$ (and $n_{j_i}<m_i$
%if $j_i=i$). 
In this case associate to $c_i$ the critical relation $(i,j_i,m_i,n_{j_i})$.

(b) $f^{m_i}(c_i)=c_j$  with $1\le j\le \nup$ and in this case associate to $c_i$ the critical relation $(i,j,n_i,0)$.

These choices ensure that properties (1) and (2) in the above definition hold. To take care that properties (3)-(6) hold
we also make the following requirement:

If both (a) and (b) hold, then only assign to $c_i$ the critical relation as in (a). If (a) holds for several $j_i\le i$, then choose
the smallest possible $j_i$ with $n_{j_i}=1$  and if there is no $j_i$ with $n_{j_i}=1$ then simply choose the smallest possible $j_i$.
Assign to $i$ only  the corresponding critical relation. Once we have done this for $i$
then repeat this construction for $i+1$.

In this way we define no new critical relation for each $1\le i \le \nup$ whose orbit is infinite and disjoint from forward orbits of $c_1,\dots,c_{i-1}$ and
from $c_1,\dots,c_\nup$, but
a unique critical relation for each of the other $i$'s.  Thus we get $N=\nup-\zeta(f)$ critical relations.

The resulting set of critical relations is realised by $f$. By construction $\mathcal F$ is proper. \end{proof}
%It is easy to see from the definition of the coordinate system $\textbf{x}$ in $\Lambda_f$ the following
%\begin{lemma}\label{nocycle}
%Assume that a finite set of critical relations $(c_{i_k}, c_{j_k}, 1, 1)$, $k=1,2,\ldots, N$, which satisfies the no-cycle condition, realized by $f$. Then the Jacobian of the map
%\begin{equation}\label{eqn:crmulti}
%g\in \Lambda_f\mapsto \{g(c_{i_k}(g))-g(c_{j_k}(g))\}_{k=1}^N
%\end{equation}
%at $g=f$ has rank $N$,
%\end{lemma}

%Given $f$, let $\zeta=\zeta(f)$ be the maximal number of critical points with pairwise disjoint infinite orbit.
\begin{proof}[Proof of Theorem~\ref{thm:crrational-minimal}]
By Remark~\ref{sigma}, it is enough to consider the case that $\sigma$ is the identity and the critical orbits avoid $\infty$.
Let $\mathcal{F}$ be a proper collection of critical relations realized by $f$.
Note that if $m,n\ge 1$ then $Q^f_{i,k;m,n}(z)$ is equal to
\begin{equation}Q^f_{i, j; m,n}(z)=\sum_{r=1}^{m-1} \frac{Df^{m-r}(f^r(c_i))}{z-f^r(c_i)}-
\sum_{s=1}^{n-1} \frac{Df^{n-s}(f^s(c_j))}{z-f^s(c_j)}\label{eq:ijmn-reduced}\end{equation}
and if $m=n=1$ then $Q^f_{i, j; m,n}(z)=0$.
By property (2) of Definition~\ref{def:lemcomb}, $(i,j;m,n)\in \mathcal{F}$ and $m,n\ge 1$ imply
\begin{equation}Df^{m-1}(f(c_i))\ne 0, Df^{n-1}(f(c_j))\ne 0\label{eq:ijmn-iter0}\end{equation}  and if $i\ne j$ then
\begin{equation} f(c_i),\dots,f^{m-1}(c_i),f(c_j),\dots,f^{n-1}(c_j)\label{eq:ijmn-iter}\end{equation} are all distinct and distinct
from $c_1,\dots,c_\nu$. Similarly, if $i=j$ then by properties (1), (2) of Definition~\ref{def:lemcomb}, $m>n$ and
$f(c_i),\dots,f^{m-1}(c_i)$, $ c_1,\dots,c_{\nu}$  are all distinct.
Hence, if $m,n\ge 1$ and $i\ne j$ then
$Q^f_{i, j; m,n}(z)$ has a non-removable pole in each of the points from the collection (\ref{eq:ijmn-iter}) and nowhere else.
In particular,   $c_1,\dots,c_\nu$ is not a pole for any $Q_{i,j;m,n}(z)$ when $m,n\ge 1$ (this holds even when $i=j$).
On the other hand, $Q^f_{i,j;m,0}(z)$ does have a pole at $c_j$ and only critical relations of this form
in $\mathcal{F}$ have a pole at $c_j$.

Suppose that the Jacobian does not have full rank. By Theorem~\ref{thm:crrational} this implies
\begin{equation}
\sum_{\substack{1\le k\le N\\ (m_k,n_k)\not=(1,1)}} a_k Q^f_{{i_k}, {j_k}; m_k, n_k}(z)=0.\label{eq:zero}\end{equation}
Let $\mathcal{F}_0$ be the set of relations $(i_k,j_k;m_k,n_k)$  in $\mathcal{F}$ in this sum for which $a_k\ne 0$
and with $(m_k,n_k)\ne (1,1)$. So (\ref{eq:zero}) is equal to the sum over the set ${\mathcal F}_0$.
By Theorem~\ref{thm:crrational}, $\mathcal{F}_0$ consists of at least one critical relation,
and obviously the properties stated in Definition~\ref{def:lemcomb}
are also satisfied for $\mathcal{F}_0$.
%By property (1)  one of the following holds:

Suppose first that there exists a critical relation $(i,j;m,0)\in \mathcal{F}_0$. In this case  by property (4) in
Definition~\ref{def:lemcomb}
there exists no $i'\ne i$, $m'>0$ so that $(i',j;m',0)\in \mathcal{F}_0$. It follows from this that $(i,j;m,0)$ is the
only term in the sum (\ref{eq:zero}) which leads to a pole at $z=c_j$. So the corresponding coefficient $a_k=0$,
a contradiction.

From now on, let us assume that for any $(i,j;m,n)\in \mathcal{F}_0$, $n>0$. Then by property (1) of Definition~\ref{def:lemcomb}, we have $i\ge j$.
Because of property (3) of Definition~\ref{def:lemcomb} we can rearrange, if necessary, the critical relations in $\mathcal{F}_0$ so that they are of the form
$(i_k,j_k;m_k,n_k)$, $1\le k\le N_0$, with $i_1<i_2<\cdots< i_{N_0}$.
If $m_k=1$ holds for all $1\le k\le N_0$, then by property (5), $(j_k, n_k)$ are pairwise distinct. Since
$Q^f_{i_k,j_k;1, n_k}$ has poles precisely at the points $f(c_{j_k}), f^2(c_{j_k}), \cdots, f^{n_k-1}(c_{j_k})$, $\sum_{k=1}^{N_0} a_k Q^f_{i_k,j_k;m_k, n_k}$ has a pole, a contradiction! So let us assume that there is a maximal $N_1\le N_0$ such that
$m_{N_1}\ge 2$.
By property (6) of  Definition~\ref{def:lemcomb}, for each $N_0\ge k>N_1$, either $j_k\not=i_{N_1}$, or $j_k=i_{N_1}$ and $n_k<m_{N_1}$. 
Together with property (2) of
Definition~\ref{def:lemcomb}, this implies that
$Q^f_{i_{k},j_k; m_k,n_k}=Q^f_{i_k,j_k;1,n_k}$ does not have a pole at $f^{m_{N_1}-1}(c_{i_{N_1}})$. For each $k<N_1$, since $i_{N_1}>i_k\ge j_k$, by property (2) of  Definition~\ref{def:lemcomb}, $Q^f_{i_k,j_k;m_k,n_k}$ does not have a pole at $f^{m_{N_1}-1}(c_{N_1})$ either. Therefore $\sum_{k=1}^{N_0} a_k Q^f_{i_k,j_k;m_k,n_k}$ has a pole at $f^{m_{N-1}-1}(c_{N_1})$, a contradiction!
\end{proof}

\section{A proof of Theorem~\ref{thm:crrational}}
\begin{remark}\label{rem:levin}
A proof of Theorem~\ref{thm:crrational} is contained essentially in~\cite{Le3}. We only outline it here (and then present another proof).
Denote $v_j(f)=f(c_j(f))$ for $j=1,\cdots,\nup$. Conjugating $f$ by a M\"obius transformation, one can assume
that $f(\infty)=\infty$, $Df(\infty)\not=0$. We label the critical values so that for some $0\le \nup'\le \nup$ the following holds:
$v_j(f)\in \C$  for  $1\le j\le \nup'$ and $v_j(f)=\infty$ for $\nup'<j\le \nup$.
Consider a subset $\Lambda_{f,\nup}\subset \text{Rat}^{\pmb \mu}_d$ of maps $g$ such that there exists $\sigma(g),b(g)\in \C$ so that $g(z)=\sigma(g)z+b(g)+O(1/z)$ as $z\to \infty$.
By~\cite{Le}, $\Lambda_{f,\nup}$ has a structure of $\nu+2$ dimensional complex manifold and
$(\sigma(g),b(g),v_1(g),\cdots,v_{\nup'}(g),v_{\nup'+1}(g)^{-1},\cdots,v_\nup(g)^{-1})$ is a holomorphic coordinate of $g\in \Lambda_{f,\nu}$. Proposition 10 of~\cite{Le3} implies
%\begin{lemma}\label{id}
that for any $(i,j;m,n)$, if $(i,j;m,n)$ is realized by $f$ and $f^m(c_i(f)), f^n(c_j(f))\not=\infty$, then
%\begin{itemize}
%\item if $n=0$, then
%\begin{equation}\label{t0}
%Q_{i,j;m,0}(x)-\widehat{Q}_{i,j;m,0}(x)=\sum_{k=1}^\nu\frac{1}{v_k-x}\frac{\partial (f^m(c_i(f))-c_j(f))}{\partial v_k}|_{f=g},
%\end{equation}
%\item if $m,n>0$, then
\begin{equation}\label{t}
\begin{array}{rl}
Q^f_{i,j;m,n}(x)- & \widehat{Q}^f_{i,j;m,n}(x)  = \\
& \\
& \sum_{k=1}^{\nup'}\frac{1}{v_k(f)-x}\frac{\partial (g^m(c_i(g))-f^n(c_j(g)))}{\partial v_k}|_{g=f} \, \, ,
\end{array}
\end{equation}
where $\widehat{Q}^f_{i,j;m,n}(x)dx^2=f_*(Q^f_{i,j;m,n}(x)dx^2)$.
%\end{itemize}
%\end{lemma}
Now Theorem~\ref{thm:crrational} can be proved by repeating
the proof of the main result of \cite{Le3} after replacing Proposition 13 of that paper by~(\ref{t}).
Instead of going into more details we give here a direct and short proof of the theorem.

\subsection{Proof of  Theorem~\ref{thm:crrational}.}\label{subsec:proofcrrational}
%\marginpar{SvS}
%In the proof below it will be convenient to have that $\infty$ is a repelling periodic point of period $>1$.
%To obtain this, pick a repelling periodic point $p$ of period  $>1$ so that no iterate of a critical point of $f$ is equal to $p$
%and replace $f$ by $M^{-1}\circ f \circ M$ where $M$ is a  M\"obius transformation so that $M(\infty)=p$.
%For the new map $\infty$ satisfies the convenient property mentioned above.
Let us first apply Thurston's pull back argument to obtain a relation of partial derivatives of $g\mapsto g^m(c_i(g))-g^n(c_j(g))$ with the quadratic differential $Q_{i,j;m,n}(z) dz^2$.
Let $L_\infty(\C)$ denote the space of all Borel measurable functions $\mu$ with $\|\mu\|_\infty<\infty$.
% and such that $\mu$ vanishes in a neighborhood of  $\infty$ and $f(\infty)$.
Note that $f^*\mu(z)=\mu(f(z))\overline{f'(z)}/f'(z)$ also belongs to the class  $L_\infty(\C)$.
\end{remark}

%\marginpar{NEW}

\begin{lemma}\label{lem:pb2diff}  Given $\mu\in L_\infty(\C)$ which vanishes in a neighborhood of  $\infty$ and $f(\infty)$,
there exists a holomorphic family $f_t$ of rational maps of degree $d$, $t\in \D_\eps$, with $f_0=f$,
and such that the following holds: For any $(i,j;m,n)$, %\marginpar{factor}
$$-\frac{1}{\pi}\int_\C (\mu-f^*\mu) Q_{i,j;m,n} |dz|^2=\left.\frac{d(f_t^m(c_i(f_t))-f_t^n(c_j(f_t)))}{dt}\right|_{t=0}.$$
\end{lemma}

\begin{proof} Assume without loss of generality that $\|\mu\|_\infty\le 1$. Then for each $t\in \D$, there are qc maps $\varphi_t, \psi_t:\C\to\C$ with complex dilatations $t\mu$ and $tf^*\mu$ respectively
such that (see e.g.~\cite{Ahl}) %\marginpar{reference added}
\begin{itemize}
\item $\varphi_t(z)=z+o(1)$, $\psi_t(z)=z+o(1)$ as $z\to\infty$ for each $t$;
\item % \marginpar{changed 2nd $\bullet$}
$f_t$ defined by $f_t\circ \psi_t=\varphi_t\circ f$ is a family of rational maps.
\end{itemize}
%It follows from the definition of the space $\Lambda_f$ that $\varphi_t, \psi_t:\C\to\C$ can %be normalized so that
%$f_t\in \Lambda_f$.
Then $\varphi_t$ and $\psi_t$ depends on $t$ holomorphically~\cite{Ahl}  %\marginpar{reference} 
and thus $\partial f_t/\bar{\partial} t=0$ in the sense of distribution, which implies that $f_t$ depends holomorphically on $t$.
Let $$\mathcal{L}_n(z)=\frac{d f_t^n(z)}{dt}|_{t=0}$$
and $L(z)=\mathcal{L}_1(z)$. Then
\begin{equation}L(x)+ Df(z) \dfrac{d}{dt}\psi_t(z) = \dfrac{d}{dt}\phi_t(f(z))\end{equation}
and therefore
\begin{equation}L(z) + Df(z) \widehat{X}(z)=X(f(z)),\label{eqn:XhatX}
\end{equation}
where $$X(z)=-\frac{1}{\pi}\int_\C \frac{\mu(\zeta)}{\zeta-z} |d\zeta|^2, \,
\widehat{X}(z)=-\frac{1}{\pi} \int_\C \frac{f^*\mu(\zeta)}{\zeta-z} |d\zeta|^2.$$
The latter formulas come from the following fact. Let $\nu\in L_\infty(\C)$ have a compact support, $\|\nu\|_\infty\le 1$ and $h_t$ ($|t|<1$) is the (unique) qc map with complex dilatation $t\nu$ such that $h_t(z)=z+o(1)$ as $z\to \infty$. Then
\begin{equation}\label{fo}
\frac{d h_t(z)}{dt}|_{t=0}=-\frac{1}{\pi}\int_\C \frac{\nu(\zeta)}{\zeta-z} |d\zeta|^2.
\end{equation}
The formula (\ref{fo}) is well-known and follows for example from the formula (6) in the proof of Theorem 1 of~\cite{Ahl}, Ch.V
(noting that a different normalisation for $h_t$ is chosen there) or by differentiating  the formula on the 2nd line of page 25 in \cite{CG}.
%where $|d(x+iy)|^2:=dxdy$.

For any $h\in \{1,2,\ldots, \nup\}$ and non-negative integer $l$, define $$S_l(c_h)=\sum_{r=1}^l Df^{l-r}(f^r(c_h)) X(f^r(c_h))$$
and $$\widehat{S}_l(c_h)=\sum_{r=1}^l Df^{l-r}(f^r(c_h)) \widehat{X}(f^r(c_h))$$
It suffices to show that for any $l,h$ as above,
\begin{equation}\label{eqn:ShatS}
S_l(c_h)-\widehat{S}_l(c_h)=\left.\frac{df_t^l(c_h(f_t))}{dt}\right|_{t=0}-\widehat{X}(f^l(c_h)).
\end{equation}
If $l=0$, then the left hand is equal to zero, and the right hand side is also equal to zero, since $\widehat{X}(c_h)=\frac{d\psi_t(c_h)}{dt}|_{t=0}$ and
$c_h(f_t)=\psi_t(c_h).$
For $l\ge 1$, we use (\ref{eqn:XhatX}):
\begin{align*}
S_l(c_h)=&\sum_{r=1}^{l}Df^{l-r}(f^r(c_h))X(f^r(c_h))\\
=&\sum_{r=1}^l Df^{l-r}(f^r(c_h)) L(f^{r-1}(c_h)) + \sum_{r=2}^l  Df^{l-r+1}(f^{r-1}(c_h))\widehat{X}(f^{r-1}(c_h)) \\
=& \mathcal{L}_l (c_h) +\sum_{r=1}^{l-1} Df^{l-r}(f^r(c_h))\widehat{X}(f^r(c_h))\\
=& \mathcal{L}_l (c_h) +\widehat{S}_l(c_h)- \widehat{X}(f^l(c_h)).
\end{align*}
Since
$$\left.\frac{df_t^l(c_h(f_t))}{dt}\right|_{t=0}=\left.\frac{df_t^l(c_h)}{dt}\right|_{t=0},$$
equation (\ref{eqn:ShatS}) follows.
\end{proof}

\begin{proof}[Proof of Theorem~\ref{thm:crrational}] Assume that the Jacobian matrix has rank less than $N$. Then there exist complex numbers $a_1, a_2, \cdots, a_N$ such that all the partial derivatives of the map
\begin{equation}\label{eqn:ak}
g\mapsto \sum_{k=1}^N a_k \left(g^{m_k}(c_{i_k}(g))-g^{n_k}(c_{j_k}(g))\right)
\end{equation}
is equal to $0$ at $g=f$. Since $\mathcal{F}$ satisfies the non-cyclic condition (\ref{eq:min11}), by Theorem~\ref{thm:cvrational}, there exists $k$ such that $(m_k,n_k)\not=(1,1)$ and $a_k\not=0$.

Given $\mu\in L_\infty(\C)$ which vanishes in a neighbourhood of $\infty$ and $f(\infty)$, %\marginpar{NEW}
let $f_t$ be given by the previous lemma.
Then for each $k=1,2,\ldots, N$, we have  %\marginpar{factor}
$$-\frac{1}{\pi}\int_\C (\mu-f^*\mu) Q_{i_k, j_k; m_k,n_k} |dz|^2=\left.\frac{d(f_t^m(c_i(f_t))-f_t^n(c_j(f_t)))}{dt}\right|_{t=0}.$$
Thus for $q$ defined as in (\ref{eqn:mainq}), 
and $\widehat{q}(z) dz^2= f_*(q(z) dz^2)$, we have
\begin{align*}
&\int_\C \mu(\hat{q}-q) |dz|^2\\
=&\int_\C (\mu-f^*\mu) q(z) |dz|^2 \\
=& -\pi\sum_{k=1}^N a_k \left.\frac{d(f_t^m(c_i(f_t))-f_t^n(c_j(f_t)))}{dt}\right|_{t=0}=0,
\end{align*}
where the %\marginpar{factor in line above} 
last equality follows from the argument in the previous paragraph.
It follows that $\widehat{q}=q$.

Assume now that $f$ is not a flexible Latt\'es example. Let us prove that $q=0$. To this end, first assume $f(\infty)\not=\infty$. Let $\varphi_i$ be the local inverse diffeomorphic branches of $f$ near $\infty$. Then $\hat{q}=q$ implies that
$$ q(z)=\sum_{i=1}^{d} q(\varphi_i(z)) \varphi_i'(z)^2$$
holds near $\infty$.
%Note that $q$ has only a finite number of poles in $\C$,
%each of which is simple and a forward iterate of a critical point of $f$.
Since $\varphi_i(\infty)\in \C$ (and is not
equal to one of the finitely many poles of $q$) and $\varphi_i'(z)=O(1/z^2)$ as $z\to \infty$,
it follows from the displayed formula $q(z)=O(1/z^4)$ at infinity. Thus $q(z) dz^2$ is an integrable meromorphic
quadratic differential. By a well-known argument, this implies that $q(z)=0$, see for example Section 3.5 of \cite{McM} and \cite{DH2}.

If $f(\infty)=\infty$, then we can find a sequence of M\"obius transformations $\sigma_l$, $l=1,2,\ldots$, converging to the identity uniformly, such that $f_{(l)}=\sigma_l\circ f \circ \sigma_l^{-1}$ satisfies $f_{(l)}(\infty)\not=\infty$ and $\infty$ avoids the critical orbits of $f_{(l)}$. Putting $g_{(l)}=\sigma_l\circ g \circ \sigma_l^{-1}$,  by (\ref{eqn:ak}), all partial derivatives of the map
$$g\mapsto \sum_{k=1}^{N} a_k \left(\sigma_l^{-1}\left(g_{(l)}^{m_k}\left(c_{i_k}(g_{(l)})\right)\right)-\sigma_l^{-1}\left(g_{(l)}^{n_k}\left(c_{j_k}(g_{(l)})\right)\right)\right)$$ 
 are equal to zero at $g=f$, hence all partial derivatives of the map
$$g\mapsto \sum_{k=1}^{N} \frac{a_k}{\sigma_{(l)}'(f^{m_k}(c_{i_k}))} ((g^{m_k}(c_{i_k}(g))-g^{n_k}(c_{j_k}(g)))$$
are equal to zero at $g=f_{(l)}$. Since $f_{(k)}$ is not a Latt\'es example, as above we obtain that
$q_{(l)}:=\sum_{k=1}^N \frac{a_k}{\sigma_{(l)}'(f^{m_k}(c_{i_k}(f)))} Q_{i_k,j_k; m_k, n_k}^{f_{(l)}}\equiv 0$.
%\marginpar{CHECK}
%where $Q_{i_k,j_k;m_k, n_k}^{f_{(l)}}$ is the expression
%defined in Subsection~\ref{subsec:32} with poles at iterates of $\sigma_k(c_1),\dots,\sigma_k(c_\nu)$.
%as above we obtain that $q_{(k)}=0$.
By continuity we conclude that $q=0$.
\end{proof}

\section{The polynomial case}\label{sec:poly}
The previous theorems also hold in the space of polynomials of degree $d$. In that case, let  $\pmb \mu=(\mu_1,\dots,\mu_\nu)$ so that
$\sum_{i=1}^\nu \mu_i= d-1$ and  let $\textbf{Pol}_d^{\pmb \mu}$ be the set of maps with critical points $c_1,\dots,c_\nu\in \C$
of orders $\mu_1,\dots,\mu_\nu$. The space $\textbf{Pol}_d^{\pmb\mu}$ is clearly an embedded submanifold of $\textbf{Rat}_d^{\hat{\pmb\mu}}$ of codimension one, where $\hat{\pmb\mu}=(\mu_1, \mu_2, \cdots,\mu_\nu, d-1)$.
%Let $\zeta_\C(f)=\zeta(f)-1$.%It follows from Theorem~\ref{thm:cvrational} that the space $\textbf{Pol}_d^{\pmb \mu}$ is a manifold,
%see Remark~\ref{rem:thm21}. The following transversality result holds:

\begin{theorem}\label{thm:poly}
Assume $f\in \textbf{Pol}_d^{\, \pmb \mu}$. Then  there exists a set
$$\mathcal F=\{(i_k, j_k; m_k, n_k), k=1,\dots,N\}\mbox{ with }N=\nu-\zeta(f)$$
of critical relations $f^{m_k}(c_{i_k})=f^{n_k}(c_{j_k})$ which are realised by $f$,
such that the Jacobian of the map
\begin{equation} \textbf{Pol}_d^{\, \pmb  \mu}\ni  g\mapsto (g^{m_k}(c_{i_k}(g))-g^{n_k}(c_{j_k}(g)))_{k=1}^N
\label{eqn:crmultipoly}\end{equation}
at $g=f$ has rank $N$.
%the Jacobian of the map
%$g\mapsto (g^{m_k}(c_{i_k}(g))-g^{n_k}(c_{j_k}(g)))_{k=1}^N$
%at $g=f$ has rank $N$.
%Assume $f\in \textbf{Rat}_d^{\pmb \mu}$ is not a flexible Latt\`es map and the critical orbits avoid the point at $\infty$. Then  there exists a set
%$$\mathcal F=\{(i_k, j_k; m_k, n_k), k=1,\dots,N\}\mbox{ with }N=\nu-\zeta(f)$$
%of critical relations $f^{m_k}(c_{i_k})=f^{n_k}(c_{j_k})$ which are realised by $f$, so that the %Jacobian of the map $\mathcal{R}_{\mathcal{F}}$ at $g=f$ has rank $N$.
%\begin{equation}\label{eqn:crmulti}
% \textbf{Rat}^{\pmb \mu}_d \ni g \mapsto \{g^{m_k}(c_{i_k}(g))-g^{n_k}(c_{j_k}(g))\}_{k=1}^N
%\end{equation}
%has rank $N$ for all $g$ close to $f$.
\end{theorem}
\begin{proof} Let $c_{\nu+1}=\infty$. For maps $g$ in $\textbf{Rat}_d^{\hat{\pmb\mu}}$ close to $f$, let $c_j(g)$ denote the critical point of $g$ close to $c_i$, $1\le i\le \nu+1$.
By Theorem~\ref{thm:crrational-minimal}, there is a set $\widehat{\mathcal{F}}=\{(i_k, j_k;m_k, n_k)\}_{j=1}^{N+1}$ of critical relations of $f$ so that the Jacobian of the map
$$\mathcal{R}_{\widehat{\mathcal{F}}}^\sigma: \textbf{Rat}_d^{\hat{\pmb \mu}}\ni g\mapsto (\sigma(g^{m_k}(c_{i_k}(g)))-\sigma(g^{n_k}(c_{j_k}(g))))_{k=1}^{N+1}$$ 
 has rank $N+1$ at $g=f$, where $\sigma$ is a M\"obius tansformation such that $\sigma(f^{m_k}(c_{i_k}))\not=\infty$ for
all $k$.
Since $\widehat{\mathcal{F}}$ is full, there is $k_0$ such that $(i_{k_0}, j_{k_0}; m_{k_0}, n_{k_0})=(\nu+1,\nu+1;1,0)$ (or $(\nu+1, \nu+1; 0,1)$). Assume without loss of generality $k_0=N+1$. Let $\mathcal{F}=\{(i_k,j_k;m_k, n_k): 1\le k\le \nu\}$ and let $\mathcal{R}$ denote the map defined by (\ref{eqn:crmultipoly}). Note that the kernel of $D_f\mathcal{R}$ is contained in the kernel of $D_f\mathcal{R}^\sigma_{\widehat{\mathcal{F}}}$, so its dimension is at most $\dim (\textbf{Rat}_d^{\hat{\pmb\mu}})-(N+1)=\dim (\textbf{Pol}_d^{\pmb\mu})-N$.  Thus the rank of
$D_f\mathcal{R}$ is at least $N$. The rank is not more than $N$, so it is equal to $N$.
% =\dim (\textbf{Pol}_d^{\pmb\mu}-$$
%In this case, the proof is almost the same as that of Theorem~\ref{thm:crrational-minimal} except that
%then we only consider critical relations associated to critical points $c_i\in \C$. The simplifying assumption
%from Remark~\ref{sigma}  that iterates of critical points avoid $\infty$ is then trivially satisfied.
%The holomorphic family $f_t$ from Proposition~\ref{lem:pb2diff} will consist of polynomials.
%As in the proof of Theorem~\ref{thm:crrational} we obtain that $\hat q=q$, and since $f(\infty)=\infty$,
%as before, we M\"obius conjugate $f$ to a rational map and conclude the proof as in the final two paragraphs of the previous section.
\end{proof}

\noindent
{\bf Acknowledgement:} The authors thank Alex Eremenko for a discussion about Theorem~\ref{thm:cvrational}, Adam Epstein
for pointing out an error in the last paragraph of section~\ref{rem:levin} in  a previous version,
%\section{A proof of Theorem~\ref{thm:crrational}} \begin{remark}\label{rem:levin}
Xavier Buff and  Lasse Rempe-Gillen for helpful comments on the introduction of the final version of this paper
and the referee for carefully reading the paper. 
This project was partly supported by the  ISF grant no: 1226/17, ERC AdG grant no: 339523 RGDD and the NSFC grant no: 11731003.

\bibliographystyle{plain}             %    other options: alpha, plain, abbrv ...

\begin{thebibliography}{1}
%\bibitem{D}
%A. Douady.
%\newblock Topological entropy of unimodal maps: monotonicity for quadratic polynomials.  {\em Real and complex dynamical systems}, 65--87
%\newblock NATO Adv. Sci. Inst. Ser. C Math. Phys. Sci., 464, Kluwer Acad. Publ., Dordrecht, 1995.
\bibitem{Ahl} L.V. Ahlfors, Lectures on quasiconformal mappings, Van Nostrand, Princeton 1966; Second Edition: A.M.S.  University Lecture Series, vol. 38, 2006.
\bibitem{Astor} M. Astorg, Summability condition and rigidity for finite type maps, Arxiv:1602.05172v1.
\bibitem{BE} X. Buff and A. Epstein, Bifurcation measure and postcritically finite rational maps,  Complex dynamics, 491--512, A K Peters, Wellesley, MA, 2009.
\bibitem{CG} L. Carleson and T.W. Gamelin, Complex dynamics, Springer Verlag, 1992.
\bibitem{DH1}  A. Douady and J.H. Hubbard, \'Etude dynamique des polyn\^omes complexes. Publications Math\'ematiques d'Orsay, 84--2.
Universit\'e de Paris-Sud, D\'epartement de Math\'ematiques, Orsay, 1984. 75 pp.
\bibitem{DH2} A. Douady and J.H. Hubbard,
A proof of Thurston's topological characterization of rational functions.  Acta Math. 171 (1993), no. 2, 263--297.
\bibitem{Ep1}  A. Epstein,  Infinitesmimal Thurston rigidity and the Fatou-Shishikura inequality, Stony Brook IMS preprint 1999\#1.
\bibitem{Ep2}  A. Epstein,  Transversality in holomorphic dynamics, \url{http://homepages.warwick.ac.uk/~mases/Transversality.pdf}.
\bibitem{Ep3} A. Epstein, Slides of talk available in
\url{https://icerm.brown.edu/materials/Slides/sp-s12-w1/Transversality_Principles_in_Holomorphic_Dynamics_\%5D_Adam_Epstein,_University_of_Warwick.pdf}
\bibitem{EG} A. Eremenko and A. Gabrielov, Rational functions with real critical points and the B. and M. Shapiro conjecture in real enumerative geometry,
Ann of Math., 155 (2002), 105--129.
\bibitem{EL} A. Eremenko and M. Lyubich, Dynamical properties of some classes of entire functions. Ann. Inst. Fourier, 42, 4 (1992), 1--32.
\bibitem{Er} A. Eremenko, A Markov-type inequality for arbitrary plane continua, Proc. AMS, 135 (2007), 1505--1510.
\bibitem{FG} C. Favre and T. Gauthier, Distribution of postcritically finite polynomials, Israel Journal of Mathematics 209 (2015), 235-292.
\bibitem{Gar} F. Gardiner, Teichmuller theory and quadratic differentials, Wiley 1987.
%and Infinitesimal Thurston Rigidity and the Fatou-Shishikura Inequality, 1999, \url{arXiv:math/9902158v1}
\bibitem{LSY} G. Levin,  M.L. Sodin and P.M.  Yuditski,  A Ruelle operator for a real Julia set. Comm. Math. Phys. 141 (1991), no. 1, 119--132.
\bibitem{Le00} G.M.Levin, On the theory of iterations of polynomial families in the complex plane. Translation from: Toeriya Funkzii, Funkzionalnyi Analiz i Ih Prilozheniya, No. 51 (1989), 94-106.
\bibitem{Le0} G. M. Levin, Polynomial Julia sets and Pade's approximations (in Russian). Proceedings of XIII Workshop on Operator's Theory in Functional Spaces (Kyubishev, 6-13 October, 1988). Kyubishev State University, Kyubishev, 1988, 113-114
\bibitem{Le1} G. Levin, On an analytic approach to the Fatou conjecture. Fund. Math. 171 (2002), no. 2, 177--196.
\bibitem{Le}  G. Levin, Multipliers of periodic orbits in spaces of rational maps. Ergodic Theory Dynam. Systems 31 (2011), 197--243.
\bibitem{Le3} G. Levin,  Perturbations of weakly expanding critical orbits. Frontiers in complex dynamics, 163--196, Princeton Math. Ser., 51, Princeton Univ. Press, Princeton, NJ, 2014.
\bibitem{LSS} G. Levin, W. Shen and S. van Strien, Monotonicity of entropy for one-parameter families of interval maps. Preprint Oct 2016.
\bibitem{Mak} P. Makienko, Remarks on the Ruelle operator and the invariant line fields problem. II.
Ergodic Theory Dynam. Systems 25 (2005), no. 5, 1561--1581.
\bibitem{McM0} C. McMullen, Amenability, Poincar\'e series and quasiconformal maps,
Invent. Math 97 (1989), 95--127.
\bibitem{McM} C. McMullen, Complex renormalisation and renormalisation, Princeton University Press, 1994, Princeton.
\bibitem{MS} W. de Melo and S. van Strien, One-dimensional dynamics,  Ergebnisse der Mathematik und ihrer Grenzgebiete. Springer-Verlag,1993,  Berlin.
\bibitem{Mil0} J. Milnor, Tsujii's monotonicity proof for real quadratic maps, unpublished 2000.
%\bibitem{MT}
%J. Milnor and W. Thurston.
%\newblock On iterated maps of the interval. {\em Dynamical systems (College Park, MD, 1986--87)},  465--563,
%\newblock Lecture Notes in Math., 1342, Springer, Berlin, 1988
\bibitem{Mil} J. Milnor, On Latt\`es maps, Dynamics on the Riemann sphere, 9--43, Eur. Math. Soc., Z\"urich, 2006.
\bibitem{MT} J. Milnor and W. Thurston,  On iterated maps of the interval. Dynamical systems (College Park, MD, 1986--87), 465--563, Lecture Notes in Math., 1342, Springer, Berlin, 1988.
\bibitem{MTr} J. Milnor and C. Tresser, On entropy and monotonicity for real cubic maps, Commun. Math. Phys. 209 (2000), 123--178.
\bibitem{Str} S. van Strien,  Misiurewicz maps unfold generically (even if they are critically non-finite). Fund. Math. 163 (2000), no. 1, 39--54.
    \bibitem{T} O. Teichm\"uller, Eine Anwendung quasikonformer Abbildungen auf das
Typenproblem, Deutsche Math., 2 (1937) 321--327. Gesammelte Abhandlungen, Springer, Berlin, 1982, 171--177.
\bibitem{Tsu0} M.  Tsujii, A note on Milnor and Thurston's monotonicity theorem. Geometry and analysis in dynamical systems (Kyoto, 1993), 60--62, Adv. Ser. Dynam. Systems, 14, World Sci. Publ., River Edge, NJ, 1994.
\bibitem{Tsu} M. Tsujii, A simple proof of monotonicity of entropy in the quadratic family, Ergodic Theory Dynam. Systems 20 (2000), 925--933.
%\bibitem{Su}
%D. Sullivan. {\em Unpublished.}
%\bibitem{T}
%M. Tsujii.
%\newblock {\em A simple proof for monotonicity of entropy in the quadratic family.}
%\newblock Ergodic Theory Dynam. Systems  {\bf 20}  (2000),  no. 3, 925--933.
\end{thebibliography}

\end{document}